\documentclass[leqno,12pt]{article}
\usepackage{amssymb,amsfonts}
\usepackage{amsmath,latexsym}

\usepackage{amstext,epic,eepic,epsf,pslatex}
\usepackage{graphicx}

\newcommand{\bepr}{{\em Proof} } 
\newcommand{\enpr}{\hfill \rule{.5em}{.5em}}

\newcommand{\R}{{\mathbb R}}
\newcommand{\N}{{\mathbb N}}

\def\XXint#1#2#3{{\setbox0=\hbox{$#1{#2#3}{\int}$ }
\vcenter{\hbox{$#2#3$ }}\kern-.6\wd0}}

\newtheorem{prop}{Proposition}[section] 
\newtheorem{thm}{Theorem}[section]

\newtheorem{cor}{Corollary}[section] 
\newtheorem{rque}{Remark}[section]

\begin{document}

\title{Source-solutions for the multi-dimensional Burgers equation}

\author{Denis Serre \\ \'Ecole Normale Sup\'erieure de Lyon\thanks{U.M.P.A., UMR CNRS--ENSL \# 5669. 46 all\'ee d'Italie, 69364 Lyon cedex 07. France. {\tt denis.serre@ens-lyon.fr}}}

\date{}

\maketitle

\begin{abstract}
We have shown in a recent collaboration that the Cauchy problem for the multi-dimen\-sional Burgers  equation is well-posed when the initial data $u(0)$ is taken in the Lebesgue space $L^1(\R^n)$, and more generally in $L^p(\R^n)$. We investigate here the situation where $u(0)$ is a bounded measure instead, focusing on the case $n=2$. This is motivated by the description of the asymptotic behaviour of solutions with integrable data, as $t\rightarrow+\infty.$
\end{abstract}



\bigskip

\paragraph{Notations.} We denote $\|\cdot\|_p$ the norm in Lebesgue $L^p(\R^n)$. The space of bounded measure over $\R^m$ is ${\cal M}(\R^m)$ and its norm is denoted $\|\cdot\|_{\cal M}$. The Dirac mass at $X\in\R^n$ is $\delta_X$ or $\delta_{x=X}$. If $\nu\in{\cal M}(\R^m)$ and $\mu\in{\cal M}(\R^q)$, then $\nu\otimes\mu$ is the measure over $\R^{m+q}$ uniquely defined by $\langle \nu\otimes\mu,\psi\rangle=\langle \nu,f\rangle\langle \mu,g\rangle$ whenever $\psi(x,y)\equiv f(x)g(y)$. The closed halves of the real line are denoted $\R_+$ and $\R_-$.

\section{Introduction}

The title of the present manuscript refers to the seminal paper \cite{LP} by Tai-Ping Liu \& Michel Pierre, where the authors studied a general one-dimensional conservation law
\begin{equation}
\label{eq:oneD}
\partial_tu+\partial_xf(u)=0,\qquad x\in\R,\,t>0
\end{equation}
when the initial data $a=u(0)$ is a bounded measure instead of a bounded or an integrable function. 

We continue here this exploration, though in a multi-dimensional context, with an equation
\begin{equation}
\label{eq:multiD}
\partial_tu+{\rm div}_x\vec f(u):=\partial_tu+\partial_1f_1(u)+\cdots+\partial_nf_n(u)=0.
\end{equation}
As in \cite{LP} we are interested, for natural reasons, in data $a$ whose total mass is finite. Whenever (\ref{eq:multiD}) is not linear, it is expected that $u$ is damped out because of dispersion. The behaviour of $\vec f(s)$ at $s=0$ is thus of great importance and we make the generic assumption (non-degeneracy) that $\vec f''(0),\ldots,\vec f^{(n+1)}(0)$ are linearly independent. Up to a change of coordinates, this amounts so saying that these vectors are parallel to those of the canonical basis. The paradigm of such conservation laws is therefore the (multi-dimensional) Burgers equation:
\begin{equation}
\label{eq:mdB}
\partial_tu+\partial_1\frac{u^2}2+\cdots+\partial_n\frac{u^{n+1}}{n+1}=0.
\end{equation}
In collaboration with L. Silvestre \cite{DSLS}, we recently proved that the Cauchy problem is well-posed in every $L^p(\R^n)$, the solution of (\ref{eq:mdB}) being instantaneously damped out as an $L^\infty$-function.

\bigskip

Following Liu \& Pierre, we are concerned with the Cauchy problem for (\ref{eq:mdB}) when the data $a$ is a bounded measure, and more precisely when $a\in{\cal M}\setminus L^1(\R^n)$. We show below that if $n\ge2$, this problem is {\em not well-posed} in the sense of Hadamard~; in particular, the Cauchy problem behaves badly at the datum $a=\delta_0$, see Corollary \ref{c:nondel}.

This negative result seems to be caused by the extreme lack of regularity of the Dirac data. When the data instead display some mild regularity in $n-1$ directions, we prove on the contrary the two following compactness results. For the sake of simplicity, we focus here on the case $n=2$. 
\begin{quote}
\begin{itemize}
\item 
Let $M,M'<\infty$ and a unit vector $\xi\in S^1$ be given. Let $K\subset L^1(\R^2)$ be the set of functions $a$ such that $\|a\|_1\le M$ and $\|\xi\cdot\nabla a\|_1\le M'$. Then the set ${\cal K}$ of solutions $u$ associated with data $a\in K$ is relatively compact in $(L^1\cap L^\infty)_{\rm loc}((0,+\infty)\times\R^2)$, and every cluster point is a true entropy solution of Burgers~; see Theorem \ref{th:pbv}.
\item Alternatively, let $a\in{\cal M}(\R^2)$ be such that
\begin{equation}
\label{eq:modc}
\lim_{h\rightarrow0}\|a(\cdot+h\xi)-a\|_{\cal M}=0.
\end{equation}
Then there exists a sequence $b_k\in L^1\cap L^\infty(\R^2)$ converging vaguely towards $a$, such that the solution $u_k(t):=S_tb_k$ converges boundedly almost everywhere over every compact subset of $(0,+\infty)\times\R^2$. Again, the limit solves the Burgers equation~; see Theorem \ref{th:pr}.
\end{itemize}
\end{quote}

Because of the density of $L^1(\R^2)$ into ${\cal M}(\R^2)$ for the vague (weak-star) topology, this compactness is expected to provide existence for some quite singular data. Because of a scaling invariance, we are especially interested in self-similar data $a=\delta_{x_1=0}\otimes g$, where $g$ is a bounded measure on the vertical axis. If uniqueness holds true, the corresponding solutions should be self-similar. We find however that severe constraints limit the set of such data for which the Cauchy problem is well-behaved. In particular, their support may not be bounded above (in the $x_2$ direction), see Proposition \ref{p:unbg}. The situation is even worse in higher space dimension.

\paragraph{Link with asymptotic behaviour.} Our study is motivated by the following observation about the two-dimensional case. Let $u$ is a solution of (\ref{eq:multiD}), with $n=2$ and $u(0)\in L^1(\R^2)$. Recall that under non-degeneracy, we may assume $f_j(s)=\frac{s^{j+1}}{j+1}+O(s^4)$. For  $\rho>\!\!>1$, the auxiliary function
$$v_\rho(t,x)=\rho u(\rho^2t,\rho x_1,x_2)$$
solves the modified equation
$$\partial_tv+\partial_1f_1^\rho(v)+\partial_2f_2^\rho(v)=0$$
where $f_j^\rho(s)=\rho^{j+1}f(v/j)=\frac{v^{j+1}}{j+1}+O(\rho^{j-3})$. The corresponding data is
$$v_\rho(0,x)=\rho a(\rho x_1,x_2),$$
 which converges in the vague topology towards a singular measure $\bar a:=\delta_{x_1=0}\otimes g$, where 
$$g(x_2)=\int_{-\infty}^{+\infty}a(x)\,dx_1.$$
When $t\rightarrow+\infty$, it is therefore tempting to compare  $u(t,x)$ to the supposed-to-be solution $\bar u$ of the Burgers equation with singular datum $\bar a$. Notice that we expect that $\bar u$ be self-similar:
$$\bar u(t,x)=t^{-1/2}V(x_1t^{-1/2},x_2)$$ where $V$ is the {\em profile}.

Because of the negative results proved below, the description of the time asymptotics turns out to be more complicated and is left aside for future work. Notice however that the situation is even more involved in space dimension $n\ge3$. For instance if $n=3$, then the scaling for the initial data is
$$a_\rho(x)=\rho a(\rho^{4/3}x_1,\rho^{1/3}x_2,\rho^{-2/3}x_3),$$
whose vague limit is just $\bar a\equiv0$.

\paragraph{Liu \& Pierre's results.}
Let us recall the main results of \cite{LP}. The initial data is always a non-negative bounded measure $\mu$, and the initial condition is interpreted in the sense that $u(t)$ converges to $\mu$ in the narrow sense as $t\rightarrow0+$. Theorem 1.1 states the uniqueness of a non-negative solution. The sign condition on the solution may be removed when the flux satifies $r\phi(r)\ge0$ (Theorem 2.2).
The existence is proved (Proposition 2.1) for a broad family of fluxes $\phi$. Finally, the link between fundamental solutions and the asymptotic behaviour is established in Theorems 3.2 and 3.3.

Notice that the uniqueness does not always hold, when we allow solutions to have a non constant sign. For instance the one-dimensional Burgers equation with a data $a=M\delta_0$ admits a one-parameter family of admissible solutions
$$u_{p,q}(t,x)=\frac xt\,{\bf 1}_{\{pt\le x\le qt\}},\qquad p\le 0\le q,\,q^2-p^2=2M,$$
among which only one has constant sign.
These so-called $N$-waves describe the asymptotic behaviour of every admissible solution of a conservation law (\ref{eq:oneD}) for which $\phi'(0)=0$ and $\phi''(0)=1$, see \cite{Daf1}.  The uniqueness can be recovered for general solutions if the initial condition is understood pointwisely (Theorem 1.3) in terms of the integrated unknown
$$v(t,x)=\int_{-\infty}^xu(t,y)\,dy.$$

Our work below reveals that there is a major gap between the one-dimensional and the two-dimensional situations. It suggests that the time-asymptotics of integrable solutions may be significantly more complex when $n\ge2$.

\paragraph{Plan of the article.}
Section \ref{s:btom} starts with well-known facts about Kru\v{z}kov's theory from \cite{Kru}. It recalls the $L^1$-theory as developped in \cite{DSLS}. Finally it describes the strategy to attack the Cauchy problem when the initial data are bounded measures.
In the short Section \ref{s:supp}, we show that the growth of the support of a solution, in the horizontal variable, is bounded by $\sqrt t\,$. Section \ref{s:mom} is two-fold. On the one hand, we prove that the sequence $u_m$ is tight, meaning that the mass cannot escape at infinity in finite time. This is done by estimating low-order moments. On the contrary, we show that high-order moments may grow arbitrarily fast, yielding to an obstruction to the Cauchy problem for arbitrary data in ${\cal M}(\R^2)$. The non-existence of self-similar solutions for many self-similar data is established in Section \ref{s:ssim}. Last but not least, the compactness of a solution set is presented in Section \ref{s:pbv}, which considers data with a directional regularity. The case of data with a modulus of continuity in one direction (\ref{eq:modc}) is treated in Section \ref{s:pr}. The Appendix gathers miscellaneous facts of smaller importance.

To prove the compactness result, we recycle a technique already used by Tartar in \cite{Tar_HW} and by Chen \& Liu in \cite{ChLi}. Since we cannot apply the {\em div-curl  Lemma}, which provided identities in terms of a Young measure associated with the sequence $u_m$, we content ourselves to use the weak upper semi-continuity result, recently established by De Rosa \& al. \cite{RST}. This furnishes inequalities only, but strong enough ones, and we are able to prove that the Young measure is almost everywhere a Dirac mass.

\bigskip

\paragraph{Acknowledgements.} The author thanks Alberto Bressan and Luis Silvestre for valuable discussions. Part of this research was done during a stay at the Department of Mathematics of Pennsylvania State University.

\section{From $L^\infty$ data to bounded measures}\label{s:btom}

\paragraph{The Kru\v{z}kov semigroup.}
Kru\v{z}kov's Theorem \cite{Kru} tells us that when $a\in L^\infty(\R^n)$, then the Cauchy problem for (\ref{eq:oneD}) admits a unique admissible solution, defined as a function $u\in L^\infty(\R_+\times\R^n)$ satisfying
\begin{equation}
\label{eq:inKru}
\int_0^\infty\int_{\R^n}(u\partial_t\phi+f(u)\cdot\nabla_x\phi)\,dxdt+\int_{\R^n}a(x)\phi(0,x)\,dx=0,\qquad\forall\phi\in{\cal D}(\R^{1+n}).
\end{equation}
as well as the so-called entropy inequalities
\begin{equation}
\label{eq:entKru}
\int_0^\infty\int_{\R^n}(|u-k|\partial_t\phi+{\rm sgn}(u-k)(f(u)-f(k))\cdot\nabla_x\phi)\,dxdt+\int_{\R^n}|a(x)-k|\phi(0,x)\,dx\ge0,
\end{equation}
for every non-negative test functions $\phi$ and every $k\in\R$. The Cauchy problem defines a semi-group over $L^\infty(\R^n)$ by $(t,a)\mapsto S_ta:=u(t,\cdot)$.

Here are the main properties of the Kru\v{z}kov's semigroup $(S_t)_{t>0}$\,:
\begin{quote}
\begin{description}
\item[Contraction.] If $a,b\in L^\infty(\R^n)$ and $b-a\in L^1(\R^n)$, then $S_tb-S_ta\in L^1(\R^n)$ and
\begin{equation}
\label{eq:contr}
\|S_tb-S_ta\|_1\le\|b-a\|_1.
\end{equation}
In particular, if $a\in L^1\cap L^\infty(\R^n)$ then $u(t)=S_ta$ remains in the space $L^1(\R^n)$ and $t\mapsto \|u(t)\|_1$ is non-increasing.
\item[Conservation of mass.] Under the same assumptions as above, we have
$$\int_{\R^n}(S_tb-S_ta)(y)\,dy=\int_{\R^n}(b-a)(y)\,dy.$$
\item[Comparison.] If $a,b\in L^\infty(\R^n)$ and $a\le b$, then $S_ta\le S_tb$.
\end{description}
\end{quote}
By using Theorem 1.2 of \cite{Cra}, and by truncating the flux, we have the important property that if $a\in L^1\cap L^\infty(\R^n)$, then $u\in C(\R_+;L^1(\R^n))$~; in particular the initial condition is satisfied in the strong sense that
\begin{equation}
\label{eq:instr}
\lim_{t\rightarrow0+}\|u(t)-a\|_1=0.
\end{equation}

\paragraph{Initial data in $L^1$.} It has been observed for a long time that the contraction property implies that the restriction of the semigroup to the subspace $L^p\cap L^\infty(\R^n)$ admits a unique continuous extension, still denoted $S_t$, to the space $L^p(\R^n)$. The case $p=1$ deserves a special attention for two reasons. On the one hand, the corresponding semigroup is contracting (see (\ref{eq:contr})). On the other hand, the $L^1$-continuity of $t\mapsto u(t)$, already known when $u(0)\in L^1\cap L^\infty(\R^n)$, extends, thanks to the contraction in the same norm, to every data in $L^1(\R^n)$. We shall therefore call $u(t):=S_ta$ the {\em abstract solution} of the Cauchy problem for (\ref{eq:multiD}) with initial data $a\in L^1(\R^n)$.

The reason why we speak of an abstract solution, instead of an admissible solution, follows from an observation made by Crandall \cite{Cra}: It is unclear whether $u$ satisfies the conservation laws (as well as the associated entropy inequalities) in the distributional sense. The knowledge that $u$ belongs to the space $C(\R_+;L^1(\R^n))$,  is not sufficient to give a meaning to the partial derivatives $\partial_jf_j(u)$~; it might happen that $f_j(u)$ is not locally integrable and thus does not defined a distribution. Of course, if the fluxes $f_j$ are globally Lipschitz, then $f\circ u$ is integrable and ones proves easily that $u$ is an admissible solution in the usual sense. This is true also if $|f(r)|=O(1+r^p)$ and $a$ is taken in $L^p(\R^m)$, because then $u\in L^\infty(\R_+;L^p(\R^n))$.

Besides these rather simple situations, it may happen that (\ref{eq:multiD}) displays a strong enough nonlinearity, which forces $u(t)$ to have a higher integrability once $t$ is positive. This is what happens in one space dimension, when the flux $f$ is a convex function, because of the dispersion relation (see \cite{Daf1}, or \cite{Daf} for an even more general statement)
$$TV_x(uf'(u)-f(u))\le\frac{2\|a\|_1}t\,,$$
which implies $u(t)\in L^\infty(\R)$ whenever $t>0$. Another useful tool is the well-known Oleinik Inequality, valid for convex fluxes,
$$f'(u)\le\frac1t\,.$$

Since these inequalities could not be generalized to several space dimensions, Crandall's concern had not been elucidated until the recent work \cite{DSLS} by L. Silvestre and the author, a paper dedicated to the multi-dimensional Burgers equation (\ref{eq:mdB}). Using the recent tool of {\em Compensated Integrability} (see \cite{Ser_DPT,Ser_CI}) and a De Giorgi-type iteration, we proved dispersion inequalities of the form
\begin{equation}
\label{eq:disp}
\|S_ta\|_q\le c_{d,p,q}t^{-\beta}\|a\|_p^\alpha,\qquad t>0,
\end{equation}
whenever $1\le p\le q\le \infty$. The exponents $\alpha,\beta$ depend upon $d,p$ and $q$ and are the only one for which the estimates are consistent with the scaling group of the equation. One deduces easily that $u(t):=S_ta$ is an entropy solution of the conservation law in the ordinary sense: one has
\begin{equation}
\label{eq:sansin}
\int_0^\infty\int_{\R^n}(u\partial_t\phi+f(u)\cdot\nabla_x\phi)\,dxdt=0,\qquad\forall\phi\in{\cal D}((0,\infty)\times\R^n),
\end{equation}
together with the corresponding entropy inequalities
\begin{equation}
\label{eq:insansin}
\int_0^\infty\int_{\R^n}(|u-k|\partial_t\phi+{\rm sgn}(u-k)(f(u)-f(k))\cdot\nabla_x\phi)\,dxdt\ge0,\quad\forall\phi\in{\cal D}_+((0,\infty)\times\R^n),\,\forall k\in\R.
\end{equation}

\begin{rque}
\label{rk:data}
We warn the reader that the dispersive estimates do not ensure that $f(u)$ be integrable in space and time up to $t=0+$. For instance, if $n=2$, one only knows
$$\int_{\R^2}|u(t,x)|^3dx=O\left(\frac1t\right).$$ 
Therefore we cannot say whether $u$ satisfies the initial condition in the integral sense of (\ref{eq:inKru}).

The function $u(t)=S_ta$ is thus a solution to the Cauchy problem in the slightly different sense that if $a\in L^p(\R^n)$, then we have $u\in C(\R_+;L^p(\R^n))$ and
\begin{equation}
\label{eq:uap}
\lim_{t\rightarrow0+}\|u(t)-a\|_p=0.
\end{equation}
To see this property in the case $p=1$, it suffices to remark that $u$ is the uniform limit in $C(\R_+;L^1(\R^n))$ of solutions $u_m=S_ta_m$ where $a_m\in L^1\cap L^\infty(\R^n)$ and $\|a_m-a\|_1\rightarrow0$. 
\end{rque}
Of course, (\ref{eq:sansin},\ref{eq:insansin}) and (\ref{eq:uap}) are sufficient to declare that $u$ is a genuine admissible solution of the Cauchy problem for the data $a\in L^p(\R^n)$.

\paragraph{Data in the space of bounded measures.}

When the initial datum $a\in{\cal M}(\R^n)$ instead, our strategy for constructing a solution of the Cauchy problem is still to approach $a$ by a sequence of data $a_m\in L^1\cap L^\infty(\R^n)$, using some regularization tool, for instance a convolution. We might also, using truncation, suppose that each $a_m$ is compactly supported. But, because $L^1(\R^n)$ is normed-closed in ${\cal M}(\R^n)$, the sequence $a_m$ is not Cauchy and thus the sequence of associated solutions $u_m(t)=S_ta_m$ is not Cauchy either. Instead, $a_m$ approximates $a$ in the mere sense of the narrow convergence, where
$$\lim_{m\rightarrow+\infty} \int_{\R^2}\psi a_mdx=\langle a,\psi\rangle,\qquad\forall \psi\in C_b(\R^n).$$
Thanks to the Uniform Boundedness Principle, $\|a_m\|_1$ is a bounded sequence. Thus the dispersive estimates show that the sequence $u_m(t):=S_ta_m$ is bounded in $L^\infty_{\rm loc}(0,\infty;L^p(\R^n))$. A natural question is thus whether $u_m$ converges in a strong enough sense that the limit solve the conservation law. A subsidiary question, though an important one, is whether this limit satisfies the initial datum in any reasonable sense.

Notice that a vague limit $u$ is actually a weak-star limit in the space $L^\infty_{\rm loc}(0,+\infty;L^p(\R^n))$, which satisfies the same estimates as (\ref{eq:disp}), that is
$$\|u(t)\|_q\le c_{d,p,q}t^{-\beta}\|a\|_p^\alpha,\qquad t>0.$$
In particular, this limit is a measurable function and the singular part of the initial data is instantaneously damped out.

\section{Width of the support of the solution}\label{s:supp}

We consider from now on the $2$-dimensional Burgers equation
\begin{equation}
\label{eq:deDB}
\partial_tu+\partial_1\frac{u^2}2+\partial_2\frac{u^3}3=0.
\end{equation}
We shall make use of the decay of the $L^1$-norm, $\|S_ta\|_1\le\|a\|_1$, together with the dispersive estimate, for which we refer to \cite{DSLS}
\begin{equation}
\label{eq:disde}
\|S_ta\|_\infty\le c_\infty t^{-\frac12}\|a\|_1^{\frac14}.
\end{equation}
\begin{prop}
\label{p:width}
Let $a\in L^1\cap L^\infty(\R^2)$ be a given initial data and $u$ be the associated admissible solution of (\ref{eq:deDB}). Suppose 
$${\rm Supp}\,a\subset[X,X']\times\R.$$
Then
\begin{equation}
\label{eq:sqrtt}
{\rm Supp}\,u(t)\subset[-c_\infty M^{\frac14}\sqrt t\,+X,X'+c_\infty M^{\frac14}\sqrt t\,]\times\R,
\end{equation}
where $M:=\|a\|_1$.
\end{prop}

\bigskip

\bepr

Without loss of generality, we may assume that $a$ is compactly supported, with support contained in the left half-plane $x_1\le0$. Let us integrate the entropy inequality
$$\partial_t|u|+\partial_1\frac{u|u|}2+\partial_2\frac{|u|^3}3\le0$$
over the domain
\begin{equation}
\label{eq:domnul}
0<t<T,\qquad x_1>c_\infty M^{\frac14}\sqrt t\,.
\end{equation}
We obtain
$$\int\!\int_{x_1>c_\infty M^{\frac14}\sqrt T}|u(T,x)|\,dx\le-\int\!\int_{\rm lateral}\left(n_t+n_1\frac u2\right)\,|u|\,ds$$
where the right-hand side is the integral over the lateral boundary $x_1=c_\infty M^{\frac14}\sqrt t\,$, and the unit normal is outward. Up to a positive factor, the quantity in parentheses equals
$$-u+c_\infty\frac{M^{\frac14}}{\sqrt t\,}\,,$$
a non-negative quantity. We deduce 
$$\int\!\int_{x_1>c_\infty M^{\frac14}\sqrt T}|u(T,x)|\,dx\le0$$
and therefore $u\equiv0$ in the domain defined by (\ref{eq:domnul}).

\enpr

\bigskip

\begin{cor}
\label{c:supp}
Let $a_m\in L^1\cap L^\infty(\R^2)$ be a sequence such that
$$\|a_m\|_1\le M\qquad\hbox{and}\qquad{\rm Supp}\,a_m\subset[X_m,X_m']\times\R$$
for some $M$ independent of $m$ and $\lim X_m=X$, $\lim X_m'=X'$. Then any weak limit $u$ of $u_m$ satisfies the same conclusion (\ref{eq:sqrtt}).
\end{cor}

\bigskip

The same technique as above can be used to prove other results of propagation. We shall use the following one.
\begin{prop}
\label{p:prop}
Let $a\in L^1\cap L^\infty(\R^2)$ be a given initial data and $u$ be the associated admissible solution of (\ref{eq:deDB}). 

If ${\rm Supp}\,a$ is contained in $\R\times\R_+$, then 
${\rm Supp}\,S_ta$ is contained in the same half-space.

If ${\rm Supp}\,a$ is contained in $\R_+\times\R$ and $a\ge0$, then 
${\rm Supp}\,S_ta$ is contained in the same half-space.
\end{prop}

\section{Moment estimates}\label{s:mom}

We recall the result established in \cite{RST}, which is valid in arbitrary space dimension $n$.
\begin{prop}
\label{p:tight}
Let $a\in L^1\cap L^\infty(\R^n)$ be given, with compact support. Then for $q\in(1,\frac{3+n}2+\frac1n)$, the functional
$$I_q[z]:=\int_{\R^n}\sum_{j=1}^n|x_j|^{\frac{q-1}j}|z(x)|\,dx$$
satisfies
$$I_q(S_ta)\le e^{c_{n,q}t}(I_q[a]+c_{n,q}t^s),\qquad\forall t>0,$$
where
$$s=1-\frac{2n(q-1)}{2+n+n^2}>0.$$
Hereabove $c_{n,q}<\infty$ is a universal constant.
\end{prop}

The role of Proposition \ref{p:tight} is to ensure the tightness of a sequence $u_m$ when the initial data $a_m$ satisfy 
$$\sup_mI_q[a_m]<\infty.$$
In other words, this prevents the mass to escape at infinity in finite time. The condition above is satisfied for instance if $a\in{\cal M}(\R^n)$ has compact support. 

The result above concerns low-order moments of the solution, in the sense that the exponents $\frac{q-1}j$ remain bounded by $\frac{1+n}2+\frac1n\,$. The exponent of $x_n$ is actually bounded by $\frac{1+n}{2n}+\frac1{n^2}\,$~; if $n=2$, this means an exponent $<1$. This contrast with the situation of higher-order moments:
\begin{prop}
\label{p:highmom}
Let $a\in L^1\cap L^\infty(\R^n)$ be given, with compact support, say that ${\rm Supp}\,a$ is contained in $[0,X]\times\R_+$. We assume $a\ge0$.

For every $\alpha>3$, the solution of the Cauchy problem to (\ref{eq:deDB}) satisfies
$$\int_{\R^2}(1+x_2)^\alpha u(t,x)\,dx\ge M+\frac{\alpha(\alpha-3)M^{\frac52}}{12c_\infty^2}\,\log\left(1+\frac{c_\infty^2t\sqrt M}{X^2}\right).$$
\end{prop}
A remarkable consequence of Proposition \ref{p:highmom} is that if a singular data $a$, concentrated along the vertical axis, is approximated by $a_m\in L^1\cap L^\infty(\R^2)$ in such a way that the horizontal width $X_m$ of the support of $a_m$ shrinks to $0$, then 
$$\lim_{m\rightarrow+\infty}\int_{\R^2}(1+x_2)^\alpha u_m(t,x)\,dx=+\infty$$
at positive times. Therefore any vague or pointwise limit $u$ satisfies
\begin{equation}
\label{eq:mominf}
\int_{\R^2}(1+x_2)^\alpha u(t,x)\,dx=+\infty,\qquad\forall t>0.
\end{equation}

\bigskip

\bepr

The weight $(1+x_2)^\alpha$ with $\alpha>3$ is a particular case of a function $\theta(x_2)$ satisfying
$$\theta>0,\qquad\theta'>0,\qquad\int_0^\infty\frac{ds}{\sqrt{\theta'(s)\,}\,}<\infty.$$
The solution is non-negative. The support of $u(t)$ is contained in $\left[0,X+c_\infty M^{\frac14}\sqrt t\,\right]$. One has
$$\frac{d}{dt}\,\int_{\R^2}\theta(x_2)\,u(t,x)\,dx=\int_{\R^2}\theta'(x_2)\,\frac{u^3}3\,dx.$$
The H\"older inequality gives us
\begin{eqnarray*}
M^3 & \le & \int_{\R^2}\theta'(x_2)\,\frac{u^3}3\,dx\cdot\left(\int_{{\rm Supp}\,u(t)}\frac{dx}{\sqrt{\theta'(x_2)}}\right)^2 \\
& \le & \left(X+c_\infty M^{\frac14}\sqrt t\right)^2\int_{\R^2}\theta'(x_2)\,\frac{u^3}3\,dx\cdot\left(\int_0^\infty\frac{ds}{\sqrt{\theta'(s)}}\right)^2 \\
& \le & 2\left(X^2+c_\infty^2 t \sqrt M\right)\int_{\R^2}\theta'(x_2)\,\frac{u^3}3\,dx\cdot\left(\int_0^\infty\frac{ds}{\sqrt{\theta'(s)}}\right)^2.
\end{eqnarray*}
This yields a differential inequality
$$\frac{d}{dt}\,\int_{\R^2}\theta(x_2)\,u(t,x)\,dx\ge\frac{\rm cst}{X^2+c_\infty^2 t \sqrt M}\,,$$
from which we derive the lower bound of the Proposition.

\enpr

\section{Scaling and self-similarity}\label{s:ssim}

The $2$-D Burgers equation (\ref{eq:deDB}) admits a scaling group
$$(t,x_1,x_2,u)\longmapsto(\mu t,\sqrt\mu\,x_1,x_2,\sqrt\mu\,u),\qquad\mu>0.$$
This means that if $u$ is an entropy solution, then
$$u^\mu(t,x):=\sqrt\mu\,u(\mu t,\sqrt\mu\,x_1,x_2)$$
is an entropy solution too. We notice that the transformation $u\mapsto u^\mu$ preserves the total mass:
$$\int_{\R^2}u^\mu(t,x)\,dx=\int_{\R^2}u(\mu t,y)\,dy=\int_{\R^2}u(0,y)\,dy.$$
It is therefore meaningful to consider self-similar solutions, which are defined as those for which the transformations above act trivially: $u^\mu\equiv u$.
Such solutions are given, in terms of their profile $W(x)=u(1,x)$, by the formula
$$u(t,x)=\frac1{\sqrt t}\,W\left(\frac{x_1}{\sqrt t}\,,x_2\right).$$
The differential equation satisfied by $W$ is (we use the letters $y_j$ to denote the self-similar variables)
\begin{equation}
\label{eq:sseq}
\partial_1\frac12(W^2-y_1W)+\partial_2\frac{W^3}3=0.
\end{equation}

If $W\in L^1(\R^2)$, then a self-similar solution admits a limit as $t\rightarrow0+$ in the sense of the narrow convergence:
$$\lim_{t\rightarrow0+}\int_{\R^2}u(t,x)\phi(x)\,dx=\lim_{t\rightarrow0+}\int_{\R^2}W(y)\phi(y_1\sqrt t\,,y_2)\,dy=\int_{-\infty}^{+\infty}g(x_2)\phi(0,x_2)\,dx_2,$$
where
\begin{equation}
\label{eq:gW}
g(x_2)=\int_\R W(s,x_2)\,ds.
\end{equation}
This initial value $a=u(0)$ is a singular measure, supported by the vertical axis:
\begin{equation}
\label{eq:agdel}
a=g(x_2)\delta_{x_1=0}=\delta_{x_1=0}\otimes g.
\end{equation}

\subsection{Constraints for self-similar data}
\label{ss:constr}

Conversely, an initial data of the form (\ref{eq:agdel}) is invariant under the same scaling, in the sense that $\langle a,\phi_\mu\rangle=\langle a,\phi \rangle$ for every $\phi\in C_c(\R^2)$ and $\mu>0$, where $\phi\mapsto\phi_\mu$ is the adjoint transformation
$$\phi_\mu(x)=\phi\left(\frac{x_1}{\sqrt\mu}\,,x_2\right).$$
If the corresponding Cauchy problem admits a unique admissible solution $u$, it must therefore be self-similar.
Anticipating on Section \ref{s:pbv}, let us make the natural assumption that $u$ is the limit of a sequence $u_m$ of solutions associated with approximate data $a_m$, whose support is contained in $[-\frac1m\,\frac1m]\times\R$, and such that $\|a_m\|_1\le\|a\|_{\cal M}=\|g\|_{\cal M}=:M$. Then Corollary \ref{c:supp} tells us that the support of $u(t)$ is contained in the strip
$$|x_1|\le c_\infty M^{\frac14}\sqrt t\,.$$
In other words, the support of the profile $W$ is contained in the strip
$$|y_1|\le c_\infty M^{\frac14}.$$
Remark that since $\|u_m(t)\|_1\le M$, we also have $\|W\|_1\le M$ and, by dispersion, $\|W\|_\infty\le c_\infty M^{\frac14}$. With H\"older inequality, we infer $\|W\|_p\le c_\infty^{1-1/p} M^{\frac14(1+3/p)}$. We can therefore estimate $g$ in various $L^p$-norms~; we start with
$$|g(x_2)|\le\int_{-c_\infty M^{\frac14}}^{c_\infty M^{\frac14}} |W(y_1,x_2)|\,dy_1\le (2c_\infty M^{\frac14})^{1-\frac1p}\left(\int |W|^pdy_1\right)^{\frac1p}.$$
We deduce
$$\|g\|_p\le (2c_\infty^2)^{1/p'} M^{\frac12(1+\frac1p)}.$$
We summarize our analysis as follows.
\begin{prop}
\label{p:contrg}
Let $a$ be a singular data of the form (\ref{eq:agdel}) with $g\in{\cal M}(\R)$. Suppose that the Cauchy problem admits a unique solution (thus a self-similar one), which is the limit of solutions associated with approximate data. Then the profile $g$ is actually a function, which satisfies
\begin{equation}
\label{eq:contrg}
\|g\|_p\le(2c_\infty^2)^{1/p'} \|g\|_1^{\frac12(1+\frac1p)},\qquad\forall p>1.
\end{equation}
\end{prop}
As a consequence, we have the following astonishing result.
\begin{cor}
\label{c:nondel}
The Cauchy problem with the datum $\delta_0$ is ill-posed in one sense or another: Either an admissible solution does not exist, or it is not self-similar (and thus not unique), or it is not the pointwise limit of entropy solutions associated with approximate data.
\end{cor}

\bigskip

More generally, we see that in order that the Cauchy problem admits a reasonable (in the sense of Proposition \ref{p:contrg}) self-similar solution, it is necessary that 
$$g\in\bigcap_{1\le p\le\infty}L^p(\R).$$
A stronger criterion can be established by examining the high-order moments. Let us suppose that $g\ge0$ with support contained in $\R_+$, so that we may choose signed approximate data: $a_m\ge0$ with support in $[-\frac1m\,,\frac1m]\times\R_+$. This implies $u_m\ge0$. If $u$ is as in Proposition \ref{p:contrg}, then $u\ge0$, that is $W\ge0$, and (\ref{eq:mominf}) tells us that
$$\int_{\R^2}(1+y_2)^\alpha W(y)\,dy=+\infty.$$
This implies that the support of $W$ is unbounded in the upper direction. With (\ref{eq:gW}) we deduce the following.
\begin{prop}
\label{p:unbg}
We make the same assumptions as in Proposition \ref{p:contrg}. In addition, we assume that $g\ge0$ and its support is bounded below (that is, $\inf\,{\rm Supp}(g)>-\infty$). Then this support is unbounded above:
$$\sup\,{\rm Supp}(g)=+\infty.$$
\end{prop}

\section{Data with partial bounded variation}\label{s:pbv}

Corollary \ref{c:nondel} suggests that we should limit ourselves to initial data that belong to some subclass of the set ${\cal M}(\R^2)$, though a  class wider  than $L^1(\R^2)$. We investigate in this section and the next one the case of initial data that display some directional regularity. We begin here with those $a\in{\cal M}(\R^2)$ that have bounded variation in some given direction.
\begin{thm}
\label{th:pbv}
Let $\xi\in S^1$ be a unit vector. Let $a\in{\cal M}(\R^2)$ be given, such that $\xi\cdot\nabla a\in{\cal M}(\R^2)$. 

Then for every sequence $a_m\in L^1\cap L^\infty(\R^2)$ which approximates $a$ in the narrow sense, and satisfies in addition
$$\sup_m\|\xi\cdot\nabla a_m\|_{\cal M}<\infty,$$
the corresponding sequence $u_m$ of solutions ($u_m(t)=S_ta_m$) is precompact in $L^p_{\rm loc}((0,+\infty)\times\R^2)$ for finite $p$. For a converging subsequence, the limit $u(t,x)$ holds boundedly almost everywhere and is an entropy solution of the Burgers equation.
\end{thm}

\paragraph{Comments.}
\begin{enumerate}
\item
The example (\ref{eq:agdel}) fulfills the assumption of the theorem provided $g$ is a function of bounded variation (choose $\xi=\vec e_2$). 
\item Theorem \ref{th:pbv} does not claim existence to the Cauchy problem with datum $a$, because it is unclear whether such a limit fits an initial condition, and if this initial value equals $a$. See our discussion in Paragraph \ref{ss:fit}.
\item Such approximate data can be taken of the form
$$a_\epsilon:=a*\rho_\epsilon,\qquad\rho_\epsilon(x)=\frac1{\epsilon^2}\rho\left(\frac x\epsilon\right),$$
which have the property that $a_{\epsilon}$ is bounded, integrable and
$$\|a_{\epsilon}\|_1\le\|a\|_{\cal M},\qquad\|\xi\cdot\nabla a_{\epsilon}\|_1=\|\rho_{\epsilon}*(\xi\cdot\nabla a)\|_1\le\|\xi\cdot\nabla a\|_{\cal M}.$$
\item 
The Theorem extends to $n\ge3$, with the assumption that $\xi_1\cdot\nabla a,\ldots,\xi_{n-1}\cdot\nabla a\in{\cal M}(\R^2)$ where $\xi_1,\ldots,\xi_{n-1}$ are linearly independent. 
\end{enumerate}

\bigskip

The rest of this section is devoted to the proof of Theorem \ref{th:pbv}.
We denote 
$$M:=\sup_m\|a_m\|_1<\infty.$$
The sequence $u_m$ is bounded in ${\cal C}(\R_+;L^1(\R^2))$, with
$$\|u_m(t)\|_1\le\|a_m\|_1\le M.$$
Because of the dispersive estimates, it is also bounded in  $L^\infty(\tau,+\infty;L^p(\R^2))$ for every $p\in[1,\infty]$ and every $\tau>0$. Up to extracting a subsequence, we may therefore assume that it converges in the Young sense: for every $g\in{\cal C}(\R)$, the sequence $(g(u_m))_{m\in\N}$ converges in the weak-star topology of $L^\infty((\tau,+\infty)\times\R^2)$ for every $\tau>0$. Denoting the limit $\bar g$, we define as usual the Young measure $(\nu_{t,x})_{(t,x)\in\R_+\times\R^2}$ by
$$\langle\nu_{t,x},g\rangle=\bar g(t,x)$$
almost everywhere. Because of the uniform bound $\|u_m\|_\infty\le c_\infty t^{-1/2}M^{1/4}$, the objects $\nu_{t,x}$ are compactly supported probabilities.

Following the strategy initiated by Tartar \cite{Tar_HW}, we shall examine the support of $\nu_{t,x}$ and prove that it is a singleton almost everywhere. This implies that $\nu_{t,x}$ is nothing but the Dirac mass at the limit $\bar u(t,x)$, and the weak-star limit commutes with continuous functions: $\bar g=g\circ \bar u.$ Then we are allowed to pass in the limit in the PDE satisfied by $u_m$, and infer that $u$ solves (\ref{eq:deDB}).

\paragraph{A tool from functional analysis.}

Unlike Tartar's one-dimensional analysis, it is not possible  to exploit the {\em div-curl lemma} of compensated compactness, because we don't have a curl-free field at our disposal. Instead, we shall apply a weak-star upper semi-continuity  result established in \cite{RST}. This result is related to the property that the map $A\mapsto(\det A)^{\frac1{d-1}}$ is {\em divergence-quasiconcave} (see \cite{Ser_DPT,Ser_CI}) over the cone ${\bf Sym}_d^+$, following a terminology coined by Fonseca \& M\"uller \cite{FM}. Let us recall that if $A:\Omega\rightarrow{\bf Sym}_d$ is a symmetric tensor defined over a $d$-dimensional domain, its row-wise divergence ${\rm Div}\,A$ is defined by
$$({\rm Div}\,A)_i=\sum_{j=1}^d\partial_ja_{ij},\qquad1\le i\le d.$$
\begin{thm}[L. De Rosa, D. S., R. Tione.]\label{th:RST}
Let $d\ge3$ be an integer and $p>d'=\frac d{d-1}$ be given. Let $\Omega$ be an open subset of $\R^d$. Let $A_\epsilon:\Omega\rightarrow{\bf Sym}_d^+$ be a sequence of positive semi-definite symmetric tensors. We assume on the one hand that $A_\epsilon\stackrel{*}\rightharpoonup A$ in $L^p(\Omega)$, and on the other hand that the sequence $({\rm Div}\,A_\epsilon)_{\epsilon>0}$ is bounded in ${\cal M}(\Omega)$.

Up to the extraction of a subsequence, let us assume that $(\det A_\epsilon)^{\frac1{d-1}}$ has a weak limit $D$ in $L^{p/d'}$. Then we have
\begin{equation}
\label{eq:usc}
D\le(\det A)^{\frac1{d-1}}.
\end{equation}
\end{thm}
Notice that if $d=2$, the div-curl Lemma tells us the stronger result that (\ref{eq:usc}) is an equality.

\bigskip

To construct such tensors, we begin by observing that 
$$\|a(\cdot+h\xi)-a\|_{\cal M}\le \|\xi\cdot\nabla a\|_{\cal M}h,\qquad\forall h\in\R.$$
By assumption, we have thus
$$\sup_m\|a_m(\cdot+h\xi)-a_m\|_1\le Ch,\qquad\forall h\in\R,$$
for some $C<\infty$. By the contraction property, we derive
$$\|u_m(t,\cdot+h\xi)-u_m(t)\|_1\le Ch,\qquad\forall h\in\R,\,\forall t>0.$$
In other words, the sequence $\xi\cdot\nabla u_m$ is bounded in $L^\infty(\R_+;{\cal M}(\R^2))$, hence in ${\cal M}((0,T)\times\R^2)$ for every $T<\infty$.

Actually, if $f:\R\rightarrow\R$ is a locally Lipschitz function, then
\begin{eqnarray*}
\|f\circ u_m(t,\cdot+h\xi)-f\circ u_m(t)\|_1 & \le & L_f(\|u_m(t)\|_\infty)\|u_m(t,\cdot+h\xi)-u_m(t)\|_1 \\ & \le & L(\|u_m(t)\|_\infty) Ch,
\end{eqnarray*}
where $L_f(s)$ is the Lipschitz constant of the restriction of $f$ to $[-s,s]$. Thanks to the dispersive estimates, this shows that the  sequence $\xi\cdot\nabla(f\circ u_m)$ is still bounded in ${\cal M}((\tau,T)\times\R^2)$ for every $0<\tau<T<\infty$.

\bigskip

For the sake of clarity, we now apply a rotation in the plane, which transforms the direction $\xi$ into $\vec e_2$. The Burgers equation (\ref{eq:mdB}) rewrites in the new coordinates
$$\partial_tu+\partial_1f(u)+\partial_2g(u)=0,$$
where $f$ and $g$ are polynomials of degree $\le3$ and valuation $\ge2$. This transformation does not alter the dispersive estimates (\ref{eq:disp}).

\bigskip

We next recall an observation made in \cite{DSLS}. Let $\eta:\R\rightarrow\R_+$ be convex (it plays the role of an entropy) with $\eta(0)=0$. Let us defined an entropy flux $\vec q$ by $q_1'(s)=f'(s)\eta'(s)$ and $q_2'(s)=g'(s)\eta'(s)$. Then the entropy production $\varepsilon_m[\eta]:=\partial_t\eta(u_m)+\partial_1q_1(u_m)+\partial_2q_2(u_m)$ is a non-negative measure, which satisfies
$$\| \varepsilon_m[\eta]\|_{{\cal M}((\tau,\infty)\times\R^2)}\le\int_{\R^2}\eta\circ u_m(\tau,x)\,dx.$$
The right-hand side is bounded in terms of the Lipschitz constant of $\eta$, by
$$L_\eta(\|u_m(\tau)\|_\infty)\|u_m(\tau)\|_1\le L_\eta(\|u_m(\tau)\|_\infty)\|a\|_1.$$
Using once more the dispersive estimates, we deduce that the sequence $\varepsilon_m[\eta]$ is bounded in ${\cal M}((\tau,\infty)\times\R^2)$ for every $\tau>0$. But since we already have a control of its last contribution $\partial_2q_2(u_m)$, we deduce that the truncated quantity $\partial_t\eta(u_m)+\partial_1q_1(u_m)$ forms a bounded sequence in ${\cal M}((\tau,T)\times\R^2)$ for every $0<\tau<T<\infty$.

By linearity, this boundedness is valid when $\eta$ is the difference of two convex functions. Because of the dispersion estimate, it even suffices that the restriction of $\eta$ to bounded intervals be the difference of two convex functions. 

\bigskip

We are now in position to define our tensors. We start with $A_m=B^+\circ u_m$, with
$$B^+(s):=\begin{pmatrix} s_+ & f(s_+) & 0 \\ f(s_+) &  q(s_+) & 0 \\ 0 & 0 & \Delta(s_+) \end{pmatrix}\in{\bf Sym}_3,$$
where $s_+=\max(s,0)$ is the positive part and
$$q(s)=\int_0^sf'(r)^2dr,\qquad\Delta(s)=sq(s)-f(s)^2.$$ 
The Cauchy--Schwarz inequality tells us that $\Delta\ge0$, and thus $B^+(s)$ is positive semi-definite.

On the one hand, we have $A_m\stackrel{*}\rightharpoonup\bar A$ in $L^\infty_{loc}$, where
$$\bar A(t,x)=\begin{pmatrix} \langle \nu_{t,x},s_+ \rangle & \langle \nu_{t,x},f(s_+) \rangle & 0 \\ \langle \nu_{t,x},f(s_+) \rangle & \langle \nu_{t,x},q(s_+) \rangle & 0 \\ 0 & 0 & \langle \nu_{t,x},\Delta(s_+) \rangle \end{pmatrix}.$$
On the other hand $\det B^+(s)=\Delta(s_+)^2$ yields
$$(\det A_m)^\frac12\stackrel{*}\rightharpoonup \langle \nu_{t,x},\Delta(s_+) \rangle.$$
Finally, we have shown above that the sequence\footnote{Here the divergence must be taken in the time and space variables.} ${\rm Div}\,A_m$ is bounded in ${\cal M}((\tau,T)\times\R^2)$ for every $0<\tau<T<\infty$.

\bigskip

Applying Theorem \ref{th:RST} with $d=3$ and $\Omega=(\tau,T)\times\R^2$, we infer that
\begin{equation}
\label{eq:young}
\langle \nu_{t,x},\Delta(s_+) \rangle^2\le\langle \nu_{t,x},\Delta(s_+) \rangle\,(\langle \nu_{t,x},s_+ \rangle\langle \nu_{t,x},q(s_+) \rangle-\langle \nu_{t,x},f(s_+) \rangle^2).
\end{equation}
Because of $\Delta\ge0$,
this yields the alternative that either $\langle \nu_{t,x},\Delta(s_+) \rangle=0$,  or
$$\langle \nu_{t,x},\Delta(s_+) \rangle\le\langle \nu_{t,x},s_+ \rangle\langle \nu_{t,x},q(s_+) \rangle-\langle \nu_{t,x},f(s_+) \rangle^2.$$
In the former case, ${\rm Supp}\,\nu_{t,x}$ is contained in the set defined by $\Delta(s_+)=0$, that is in the negative axis $\R_-$.
The latter situation instead can be recast as $\langle \nu_{t,x}\otimes \nu_{t,x},F \rangle\le0$, where
\begin{eqnarray*}
F(r,s) & := & \Delta(r_+)+\Delta(s_+)-r_+q(s_+)-s_+q(r_+)+2f(r_+)f(s_+) \\
& = & (r_+-s_+)\int_{s_+}^{r_+}f'(\sigma)^2d\sigma-\left(\int_{s_+}^{r_+}f'(\sigma)\,d\sigma)\right)^2.
\end{eqnarray*}
By Cauchy-Schwarz, $F$ is non-negative and vanishes only if $r_+=s_+$ (remark that $f'$ is not constant on an interval). Since $\nu_{t,x}\otimes \nu_{t,x}$ is a probability, our inequality implies that the support of $\nu_{t,x}\otimes \nu_{t,x}$ is contained in the planar subset defined by $F=0$, that is $r_+=s_+$.

\bigskip

The conclusion  of this calculation is that
$${\rm Supp}\,\nu_{t,x}\otimes \nu_{t,x}\subset\R_-^2\bigcup L,$$
where $L$ denotes the diagonal ($r=s$) of $\R^2$. The same calculation can be carried out, starting from the entropy $\eta(s)=s_-:=\max(-s,0)$ instead of $s_+$. 
This yields now the symmetric result that
$${\rm Supp}\,\nu_{t,x}\otimes \nu_{t,x}\subset\R_+^2\bigcup L.$$
Combining both results, we conclude that
$${\rm Supp}\,\nu_{t,x}\otimes \nu_{t,x}\subset\left(\R_-^2\bigcup L\right)\bigcap\left(\R_+^2\bigcup L\right)=L.$$
Since 
$${\rm Supp}\,\nu_{t,x}\otimes \nu_{t,x}=({\rm Supp}\,\nu_{t,x})\times ({\rm Supp}\,\nu_{t,x}),$$
this tells us that ${\rm Supp}\,\nu_{t,x}$ is a singleton, as expected.

\bigskip

This proves that our sequence $u_m$ converges pointwise almost everywhere. Because of the $L^\infty$ bound on $(\tau,\infty)\times\R^2$ for every $\tau>0$, and using dominated convergence, we may pass to the limit in the identities
$$\int_0^\infty\!\!\int_{\R^2}\left(u_m\partial_t\phi+\frac{u^2_m}2\,\partial_1\phi+\frac{u^3_m}3\,\partial_2\phi\right)\,dy\,dt=0,\qquad\forall\phi\in{\cal D}((0,\infty)\times\R^2)$$
as well as in the entropy inequalities. We infer
\begin{equation}
\label{eq:wk}
\int_0^\infty\!\!\int_{\R^2}\left(u\partial_t\phi+\frac{u^2}2\,\partial_1\phi+\frac{u^3}3\,\partial_2\phi\right)\,dy\,dt=0,\qquad\forall\phi\in{\cal D}((0,\infty)\times\R^2),
\end{equation}
together with the entropy inequalities. This ends the proof of Theorem \ref{th:pbv}.

\subsection{The problem of fitting the initial condition}
\label{ss:fit}

As mentionned in Remark \ref{rk:data}, we don't know whether $u^3$ is integrable over $(0,1)\times B$ when $B\subset \R^2$ is a bounded set. We even suspect that it is not so in reasonnable situations, namely those of a self-similar solution. This raises a difficulty regarding the initial data. On the one hand, when $a\in{\cal M}\setminus L^1(\R^2)$, $a_m$ does not converge in norm towards $a$~; we infer that $u_m$ is not a Cauchy sequence in $C(\R_+;{\cal M}(\R^2))$. On the other hand, if $\phi\in{\cal D}(\R^2)$, we may not pass to the limit in
$$\int_{\R^2}u_m(t,x)\phi(x)\,dx-\int_{\R^2}a_m(x)\phi(x)\,dx=\int_0^t\!\int_{\R^2}\left(\frac{u_m^2}2\,\partial_1\phi+\frac{u_m^3}3\,\partial_2\phi\right)\,dxds,$$
because of the lack of uniform integrability of the last term.

We thus don't know so far whether $u$ fits the initial data $a$ in any sense. We don't even exclude the possibility that $u(t)$, which is bounded in $L^1(\R^2)$, would not converge at all as $t\rightarrow0+$. In such a case, it would admit a non-unique cluster value for the narrow convergence.

Even if $\lim_{t\rightarrow0+}u(t)=:u(0)$ exists in the narrow sense, we don't know whether it equals $a$. The Section \ref{s:ssim} gives us some clues. Let us consider a datum $a$ of the self-similar form (\ref{eq:agdel}). Choosing an $a_m$ whose support shrinks to the vertical axis as $m\rightarrow+\infty$, we know (Corollary \ref{c:supp}) that the support of $u(t)$ shrinks to the axis as $t\rightarrow0+$, its width being an $O(\sqrt t\,)$. If $u(0)$ exists, it is therefore supported by the vertical axis, thus self-similar, and it is reasonnable to suppose that $u$ is self-similar too. But then $u(0)$ has to satisfy the constraints described in Paragraph \ref{ss:constr}. If $a$ does not satisfy them, then certainly $u(0)\ne a$. 

\bigskip

The only positive result that we have in this direction is the following.
\begin{prop}
\label{p:ciav}
With the same assumptions as in Theorem \ref{th:pbv}, we have
$$\lim_{t\rightarrow0+}\int_{\R^2}u(t,x)\phi(x_1)\,dx=\langle a,\phi\otimes{\bf1}\rangle.$$
\end{prop}
In other words, the integral of $u$ along vertical lines satisfies the correct relation:
$$\lim_{t\rightarrow0+}\int_\R u(t,x)\,dx_2=\pi_* a$$
where $\pi:x\mapsto x_1$ is the vertical projection and the right-hand side is the pushforward measure. The proof consists simply in passing to the limit into
$$\int_{\R^2}u_m(t,x)\phi(x_1)\,dx-\int_{\R^2}a_m(x)\phi(x_1)\,dx=\int_0^t\!\int_{\R^2}\frac{u_m^2}2\,\phi'(x_1)\,dxds,$$
where now every term is under control ; the last one because of the dispersive estimate $\|u_m(t)\|_2\le c_\infty^{1/2}t^{-1/4}M^{5/8}$.

Of course, Proposition \ref{p:ciav} is of very little interest for data of the form (\ref{eq:agdel}), since it tells only that $u(0)$ has the right mass.

\section{Data with partial regularity}\label{s:pr}

We extend now the compactness result to a more general class of initial data. 
\begin{thm}
\label{th:pr}
Let $\xi\in S^1$ be a given unit vector. Let $a\in{\cal M}(\R^2)$ be such that 
$$\lim_{h\rightarrow0}\|a(\cdot+h\xi)-a\|_{\cal M}=0.$$
Then there exists a sequence $b_k$ in $L^1\cap L^\infty(\R^2)$, converging towards $a$ in the vague topology, and such that the solution sequence ($u_k(t):=S_tb_k$) converges boundedly almost everywhere on every compact subset of $(0,+\infty)\times\R^2$.

The pointwise limit $u$ of $u_k$ is an entropy solution of (\ref{eq:deDB}), which satisfies the dispersion inequalities
$$\|u(t)\|_1\le\|a\|_{\cal M},\qquad\|u(t)\|_\infty\le c_\infty t^{-1/2}\|a\|_{\cal M}^{1/4}.$$
\end{thm}

\bigskip

Theorem \ref{th:pr} applies to the important case of a data of the form (\ref{eq:agdel}) when the measure $g$ is actually an $L^1$-function, because then
$$\|a(\cdot+h\vec e_2)-a\|_{\cal M}=\|g(\cdot+h)-g\|_1\stackrel{h\rightarrow0}\longrightarrow0.$$

\bigskip

\bepr

Let $a_\epsilon=\theta_\epsilon*_\xi a$ be an approximate data obtained by convolution in the direction $\xi$ only. The assumption ensures that on the one hand $\xi\cdot\nabla a_\epsilon\in{\cal M}(\R^2)$, and on the other hand 
\begin{equation}
\label{eq:contrM}
\lim_{\epsilon\rightarrow0+}\|a_\epsilon-a\|_{\cal M}=0.
\end{equation}
Theorem \ref{th:pbv} tells us that for every $\epsilon>0$, there exists a sequence of data $b_{\epsilon m}\in L^1\cap L^\infty(\R^2)$, which tends to $a_\epsilon$ for the narrow topology, and such that the corresponding sequence of solutions $(u_{\epsilon m})_m$ converges boundedly almost everywhere towards a solution $u_\epsilon$  of (\ref{eq:deDB}).

Recall that in the proof of Theorem \ref{th:pbv}, one could have chosen {\em a priori} a smooth kernel $\rho_k(x)=k^2\rho(kx)$ and set $a_{\epsilon k}:=\rho_k*a_\epsilon$. Then $b_{\epsilon m}$ is a subsequence, meaning that $b_{\epsilon m}=a_{\epsilon \phi(m)}$ for some increasing function $\phi$. By a diagonal selection, we may actually take the same function $\phi$ for every parameters $\epsilon=\frac1N$ for $N\in\N$.

We remark now that
$$\sup_t\|u_{\epsilon m}(t)-u_{\nu m}(t)\|_1\le\|a_{\epsilon m}-a_{\nu m}\|_1\le\|a_\epsilon-a_\nu\|_{\cal M}\stackrel{\epsilon,\nu\rightarrow0}\longrightarrow0.$$
Reminding the uniform estimate $\|u_{\epsilon m}(t)\|_\infty\le c_\infty t^{-1/2}M^{1/4}$ with an $M$ independent of $\epsilon$ and $m$ (one may take $M=\|a\|_{\cal M}$), we apply Dominated Convergence to
$$\int_K|u_{\epsilon m}-u_{\nu m}|(t,x)\,dxdt$$
for every relatively compact domain  $K\subset(0,\infty)\times\R^2$, and we obtain
$$\int_\tau^T\!\int_K|u_{\epsilon }-u_{\nu }|(t,x)\,dxdt\le(T-\tau)\|a_\epsilon-a_\nu\|_{\cal M}.$$
this in turn implies
$$\int_\tau^T\!\int_{\R^2}|u_{\epsilon }-u_{\nu }|(t,x)\,dxdt\le(T-\tau)\|a_\epsilon-a_\nu\|_{\cal M},$$
that is
$$\sup_t\|u_{\epsilon }(t)-u_{\nu }(t)\|_1\le\|a_\epsilon-a_\nu\|_{\cal M}.$$
The sequence $u_{1/N}$ is thus time-uniformly convergent in $L^1(\R^2)$. Its limit $u$ is clearly  a solution of (\ref{eq:deDB}), and the convergence is pointwise almost everywhere, up to a new extraction.

Another diagonal selection yields a converging subsequence of $u_{1/N,N}$, towards $u$.

\enpr

\section*{Appendix}

We gather miscellaneous facts of smaller importance.

\subsection*{Characteristics of self-similar flows}

The characteristics of (\ref{eq:sseq}) obey to the differential system
$$\dot x_1=W-\frac{x_1}2,\qquad \dot x_2=W^2,\qquad \dot W=W,$$
which can be integrated into
$$x_1W-x_2={\rm cst},\qquad(x_1-W)W={\rm cst.}$$
This defines a two-parameters family of curves, which can be used to build an admissible solution of (\ref{eq:sseq}). Of course, a compactly supported solution will involve shock waves, for which we have to write the Rankine--Hugoniot relation and the Oleinik entropy inequalities, both varying with the normal direction.

\subsection*{Very self-similar flows}

The equation (\ref{eq:sseq}) admits a class of local solutions which depend upon one scalar variable only, namely
$$W(y)=y_1V\left(\frac{y_2}{y_1^2}\right).$$
The reduced profile obeys now an ODE
\begin{equation}
\label{eq:vss}
(V^2+z(1-2V))V'=V(1-V),
\end{equation}
where $z$ stands for the remaining argument. Remarkably enough, this equation can be integrated explicitly, into
$$\frac{V^2-z}{V(1-V)}={\rm cst.}$$
The corresponding solution of the Burgers equation is given implicitely by
\begin{equation}
\label{eq:vsw}
W^2+c(y_1-W)W=y_2
\end{equation} 
or equivalently
$$tu^2+c(x_1-tu)u=x_2$$
for some constant $c\in\R$.

We notice that the mass of such a function is infinite unless the domain is truncated:
$$\int u(t,y)\,dy=\int W(y)\,dy=\int y_1V\left(\frac{y_2}{y_1^2}\right)\,dy=\int y_1^3V(z)\,dy_1dz.$$
A truncation must involve shock waves. Again, we remark that the shock locus separating a very self-similar flow from the rest state ($u\equiv0$) can be completely described. The Rankine--Hugoniot relation tells us that the shock locus obeys to the differential equation
$$\frac{W-y_1}2\,dy_2-\frac{W^2}3\,dy_1=0.$$
This can be integrated into
\begin{equation}
\label{eq:vsss}
W^2\left(c(W-y_1)^2+W(\frac{4y_1}3-W)\right)={\rm cst}.
\end{equation}
The shock locus can then be determined by eliminating $W$ between (\ref{eq:vsss}) and (\ref{eq:vsw}). Notice the situation encountered when $c=\frac43$\,, where (\ref{eq:vsss}) simplifies into
$$W(W-2y_1)={\rm cst.}$$

\section*{Full USC}

When we drop the regularity assumption of either Theorem \ref{th:pbv} or \ref{th:pr}, and we suppose  $a\in{\cal M}(\R^2)$ and $a\ge0$, we can only apply the upper semi-continuity to a tensor $T\circ u_m$ given by 
$$T=((g_{ij}))_{0\le i,j\le2}=(g_{0\bullet},g_{1\bullet},g_{2\bullet})$$
where each line is made of an entropy and its fluxes, that is 
$$g'(s)=h(s)X(s)\otimes X(s),\qquad X(s):=\begin{pmatrix} 1 \\ f_1'(s) \\ f_2'(s) \end{pmatrix}.$$
The positivity is ensured by $g(0)=0$ and $h>0$.

Considering a sequence of solution $u_m$ associated with approximate solutions, controlled in $L^\infty_{\rm loc}$, a Young measure satisfies
$$\langle\nu,\sqrt{\det T\,}\rangle^2\le\det\langle\nu,T\rangle.$$
By polarizing, this is equivalent to the inequality
\begin{equation}
\label{eq:inNF}
\langle N,G(u,v,w)\rangle\le0,\qquad N:=\nu\otimes\nu\otimes\nu,
\end{equation}
where the symmetric function $F$ is given by
\begin{eqnarray*}
G(u,v,w) & = & 2\left(\sqrt{\Delta(u)\Delta(v)\,}+\sqrt{\Delta(v)\Delta(w)\,}+\sqrt{\Delta(w)\Delta(u)\,}\right) \\
& & -\det(g_{0\bullet}(u),g_{1\bullet}(v),g_{2\bullet}(w))-\cdots,\qquad \Delta:=\det T.
\end{eqnarray*}
The dots in the formula above stand for five other similar terms.

As in the case studied in Section \ref{s:pbv}, we have $G\equiv0$ and $\nabla G\equiv0$ along the diagonal, here defined by $u=v=w$. It is however not as flat (remember that we had $\nabla^2F\equiv0$ and $\nabla^3F\equiv0$ along the diagonal $u=v$), since
$$G_{uu}(u,u,u)=-\frac{\Delta'(u)^2}{\Delta(u)}=:2\alpha(u)<0,\qquad G_{uv}(u,u,u)=\alpha(u).$$ 
We infer that 
$$G(u,v,w)\sim\alpha\left(\frac{u+v+w}3\right)\,(u^2+v^2+w^2-uv-vw-wu)$$
in a neighbourhood of the diagonal. Since this expression is non-positive, one cannot conclude from Inequality (\ref{eq:inNF}). For instance, a measure $\nu$ with a small enough support certainly satisfies (\ref{eq:inNF}).


\begin{thebibliography}{00}






\bibitem{ChLi} G.-Q. Chen, Y.-G. Liu. A study on application approaches of the theory of compensatedcompactness. {\em Chinese Science Bull.}, {\bf 34} (1989), pp 15--19.
%

\bibitem{Cra} M. Crandall. The semigroup approach to first order quasilinear equations in several space variables. {\em Israel J. Math.}, {\bf 12} (1972), pp 108--132.



\bibitem{Daf} C. Dafermos. Regularity and large time behaviour of solutions of a conservation law without convexity. {\em Proc. Royal Soc. Edinburgh}, {\bf 99A} (1985), pp 201--239.

\bibitem{Daf1} C. Dafermos. Asymptotic behaviour of solutions of hyperbolic balance laws. In: Bardos C., Bessis D. (eds) {\em Bifurcation Phenomena in Mathematical Physics and Related Topics}, NATO Advanced Study Institutes Series (Series C -- Mathematical and Physical Sciences), vol 54. Springer, Dordrecht (1980).




\bibitem{FM} I. Fonseca, S. M\"uller. ${\cal A}$-Quasiconvexity, lower semicontinuity, and Young measures. {\em SIAM J. Math. Anal.}, {\bf30} (1999), pp 1355--1390.


        
\bibitem{Kru} S. Kru\v{z}kov. First order quasilinear equations with several independent variables (in Russian). {\em Mat. Sbornik (N.S.)}, {\bf 81 (123)} (1970), pp 228--255.
        
\bibitem{Lax} P. Lax. Hyperbolic systems of conservation laws, II. {\em Comm. Pure Appl. Math.}, {\bf 10} (1957), pp 537--566.
        


\bibitem{LPT} P.-L. Lions, B. Perthame, E. Tadmor. A kinetic formulation of multidimensional scalar conservation laws and related equations. {\em J. Amer. Math. Soc.}, {\bf 7} (1994), pp 169--191.
        
\bibitem{LP} T.-P. Liu, M. Pierre. Source-solutions and asymptotic behavior in conservation laws. {\em J. Differential Eq.}, {\bf 51} (1984), pp 419--441.
        
        

\bibitem{RST} L. De Rosa, D. Serre, R. Tione. On the upper semicontinuity of a quasiconcave functional. {\em Submitted.} Preprint arXiv:1906.06510.
        


\bibitem{Ser_DPT} D. Serre. Divergence-free positive symmetric tensors  and fluid dynamics. {\em Annales de l'Institut Henri Poincar\'e (analyse non lin\'eaire),}  {\bf35} (2018), pp 1209--1234. 
      

\bibitem{Ser_CI} D. Serre. Compensated integrability.  Applications to the Vlasov--Poisson equation and other models in mathematical physics. {\em Journal de Math\'ematiques Pures et Appliqu\'ees,} {\bf 127} (2019), pp 67--88.

\bibitem{DSLS} D. Serre, L. Silvestre. Multi-dimensional scalar conservation laws with unbounded initial data: well-posedness and dispersive estimates. {\em Arch. Rat. Mech. Anal.}, {\bf234} (2019), pp 1391--1411.

\bibitem{Tar_HW} L. Tartar. Compensated compactness and applications to partial differential equations. {\em Nonlinear analysis and mechanics: Heriot-Watt Symposium}, Vol. IV, Res. Notes in Math., 39, Pitman (1979), pp 136--212.
        

      



 

\end{thebibliography}
\end{document}